\documentclass[a4paper]{amsart}
\usepackage{amssymb}

\title[Curves with many points]{New curves with many points over small finite fields}
\author{Karl R\"okaeus}
\address{Karl R{\"o}kaeus \\ Korteweg de Vries Instituut voor Wiskunde \\ Universiteit van
Amsterdam \\ P.O. Box 94248  \\1090 GE Amsterdam \\ The Netherlands}
\email{S.K.F.Rokaeus@uva.nl}
\date{}

\newcommand{\comment}[1]{}

\DeclareMathOperator{\PP}{\mathbf{P}}
\DeclareMathOperator{\Cl}{Cl}

\DeclareMathOperator{\NN}{N}
\DeclareMathOperator{\FF}{\mathbf{F}}
\DeclareMathOperator{\ZZ}{\mathbf{Z}}

\DeclareMathOperator{\Supp}{Supp}

\theoremstyle{remark}

\newtheorem*{example}{Example}
\newtheorem*{acknowledgements}{Acknowledgments}

\begin{document}
\maketitle

\emph{Abstract.}\footnote{The author is supported by The Wenner-Gren
Foundations Postdoctoral Grant}
We use class field theory to search for curves with many rational points over the finite fields of cardinality $\leq5$. By going through abelian covers of each curve of genus $\leq2$ over these fields we find a number of new curves. In particular, over $\FF_2$ we settle the question of how many points there can be on a curve of genus 17 by finding one with 18 points. The search is aided by computer; in some cases it is exhaustive for this type of curve of genus up to 50.

\section{Introduction}
Let $q$ be a prime power and let $g$ be a non-negative integer. Denote by $\NN_q(g)$ the maximum number of rational points possible on a projective, smooth and geometrically irreducible curve of genus $g$, defined
over the finite field $\FF_q$; equivalently, the number of rational places possible in a global function field of genus $g$ with full constant field $\FF_q$. Interest in determining $\NN_q(g)$ arose in the 1980's, among others with \cite{Serre}. It has turned out to be a difficult task; it is still open even for small genera over $\FF_2$, starting with $g=12$ and $16$. See \cite{vanderGeerVlugt} for an overview of the history of the question and of different methods that have been used to investigate it. The current intervals in which $\NN_q(g)$ is known to lie for $g\leq 50$ are collected on the webpage \cite{manyPoints}. In this article we improve upon the lower bounds of many of these intervals by constructing new examples of curves with many points.

Most of the lower bounds in \cite{manyPoints} are obtained as abelian covers of some curve of low genus; often by going through for example Artin-Schreier or Kummer extensions, or fibre products of such extensions, but also by using class groups. We utilize the ray class groups to search through abelian extensions of curves of small genus in a systematic fashion; the aim was to find all improvements of \cite{manyPoints} that can be obtained in this way. More precisely, for each finite field $k$ of cardinality $\leq 5$ and for each curve $C/k$ of genus $\leq 2$ we perform a computer search that check all abelian covers of $C$ whose conductor degree is bounded by $d$, where the limit $d$ is set in each case by our computer power.

Over $\FF_2$ our search is exhaustive for the base curves of genus 2 with $\leq 4$ points. Here it yields a curve that shows $\NN_2(17)=18$, and improves upon the lower bounds also for genus $45$, $46$ and $48$. Over $\FF_q$ for $q=3,4,5$ we also found a number of improvements of \cite{manyPoints}, see Section \ref{results}. As expected, the search also reproduced most of the old records in \cite{manyPoints}.

The search was done using Magma \cite{magma}. We have formulated everything in terms of function fields (rather than curves), see \cite{Stich} and \cite{NiederXing}.

\section{Method}
\subsection*{Notation}
By a function field $F/k$ we will always mean a global function field with full constant field $k$, i.e., $k$ is algebraically closed in $F$. We write $\PP_F$ for the set of places of $F$.
\subsection*{Class Field Theory}
Let $F/k$ be a function field, let $D$ be an effective divisor of $F$ and let $\Cl_D(F)$ (or just $\Cl_D$) be the ray divisor class group modulo $D$ of $F$, i.e., the group of divisors relatively prime to $D$ modulo principal divisors congruent to 1 modulo $D$. There is an inclusion reversing correspondence between finite abelian extensions $K/F$ with conductor $\leq D$, and subgroups $U\subset \Cl_D(F)$ of finite index, see \cite{Hess}.
For $U$ a subgroup of $\Cl_D$ of finite index $d$, let $F^D_U/F$ be the corresponding abelian extension. This extension has degree $d$; it is unramified outside of the support of $D$, where a place splits completely if and only if its image lies in $U$. When $S\subset\PP_F$ is a set of places whose images in $\Cl_D$ generate $U$ we also write $F^D_S/F$ for $F^D_U/F$. This is then the largest abelian extension of $F$ with conductor $\leq D$ such that all places in $S$ split completely.

The ramification behavior of $F^D_U/F$ at the places in $\Supp(D)$ can be computed from the images of $U$ in $\Cl_{D'}(F)$ for $D'\leq D$. Algorithms for computing the ray class groups are given in \cite{Hess}, where it is also indicated how to compute different invariants of the class field. This is implemented in Magma \cite{magma}.

\subsection*{Organization of the Search}
We work over $k=\FF_q$ for $q=2,3,4,5$. Using Magma, for all function fields $F$ of genus 1 and 2 and for various divisors $D$  we construct $\Cl_D(F)$; for each rational place not in the support of $D$ we take the quotient with the subgroup generated by it, list the subgroups of these finite quotients and pull them back to $\Cl_D(F)$. In this way we get all abelian  extensions $K/F$ of conductor $\leq D$ in which at least one rational place splits completely (so that $k$ is the full constant field of $F$). We compare this to the tables \cite{manyPoints} and record the improvements that we find.

Our computer power puts a limit on the degree of the divisors that we can check with this method (depending on the size of the constant field of $F$ and also on its number of places of different degrees). However, for some curves over $\FF_2$ it is possible to check \emph{all} abelian extensions that could give an improvement of the old records for $g\leq 50$: Let $F/\FF_2$ be a function field of genus 2 with $N$ rational places. For triples $(d,s,m)$ of non-negative integers where $0<s\leq N$ and $0\leq m\leq N-s$, consider an abelian extension $K/F$ of degree $d$, corresponding to a subgroup of $\Cl_D(F)$ where $D$ has $m$ rational places in its support, such that $s$ rational places of $F$ splits completely in $K$. Then there are at least $ds$ rational places in $K$, and so we only consider $d$ such that $ds\leq 40$ as this is the maximum number of points on a curve of genus $\leq 50$ over $\FF_2$. Moreover, the datum $(d,s,m)$ gives an upper bound on the number of rational places in $K$; except for the $ds$ places that split completely we might get a contribution from places in the support of $D$ as big as $d_{max}m$ where $1\leq d_{max}<d$ is the greatest proper divisor of $d$. Comparing this number to the tables \cite{manyPoints} to see for which genera it would give a new record, we get an upper bound $G$ for the genus of $K$. The last formula in \cite{Hess} now gives an upper bound $B$ on the degree of $D$ (which in case $d$ is a prime is just $\deg(D)\leq B=2(G-1-d)/(d-1)$). So for the triple $(d,s,m)$, for each divisor $D$ of degree $\leq B$ with $m$ rational places in its support we construct its ray class group $\Cl_D(F)$. Then, for each set $S$ of cardinality $s$ consisting of rational places of $F$ not in the support of $D$, we check all extensions that correspond to a subgroup of $\Cl_D$ that contains the subgroup generated by the places in $S$ and is of index $d$. For $N\leq 4$ the bounds on $\deg(D)$ are good enough that we are able to check all possible extensions of these genus 2 fields (leaving the possibility that there are new curves among the abelian extensions of the three genus 2 fields with 5 rational places or the one with 6 rational places).

\section{Results}\label{results}
In this section we state the new records that we found, together with the intervals in which $\NN_q(g)$ was known to lie prior to their submission to \cite{manyPoints}. These new curves, together with an earlier version of this report, have appeared in \cite{manyPoints}. I was later informed by Ducet that he and Fieker independently are performing a similar search, which also yields a curve of genus 17 with 18 points.

The details required to construct these curves are given in Section \ref{details}.

\noindent\begin{minipage}{\textwidth}
\vspace{.2in}\noindent{\large
New entries for the tables over $\mathbf{F}_2$.}\\[.05in]
\smallskip
\vbox{
\bigskip\centerline{\def\quad{\hskip 0.6em\relax}
\def\quod{\hskip 0.5em\relax }
\vbox{\offinterlineskip
\hrule
\halign{&\vrule#&\strut\quod\hfil#\quad\cr
height2pt&\omit&&\omit&&\omit&\cr
& $g$ && $N$ &&  Interval  &\cr
\noalign{\hrule}
& $17$ && $18$ && $[17,18]$  &\cr
& $45$ && $36$ && $[33,37]$  &\cr
& $46$ && $36$ && $[34,38]$  &\cr
& $48$ && $35$ && $[34,39]$  &\cr
} \hrule}
}}
\end{minipage}

\noindent\begin{minipage}{\textwidth}
\vspace{.2in}\noindent{\large
New entries for the tables over $\mathbf{F}_3$.}\\[.05in]
\smallskip
\vbox{
\bigskip\centerline{\def\quad{\hskip 0.6em\relax}
\def\quod{\hskip 0.5em\relax }
\vbox{\offinterlineskip
\hrule
\halign{&\vrule#&\strut\quod\hfil#\quad\cr
height2pt&\omit&&\omit&&\omit&\cr
& $g$ && $N$ &&  Interval  &\cr
\noalign{\hrule}
& $17$ && $28$ && $[25,30]$  &\cr
& $22$ && $33$ && $[30,36]$  &\cr
& $33$ && $48$ && $[46,49]$  &\cr
& $46$ && $60$ && $[55,63]$  &\cr
} \hrule}
}}
\end{minipage}

\noindent\begin{minipage}{\textwidth}
\vspace{.2in}\noindent{\large
New entries for the tables over $\mathbf{F}_4$.}\\[.05in]
\smallskip
\vbox{
\bigskip\centerline{\def\quad{\hskip 0.6em\relax}
\def\quod{\hskip 0.5em\relax }
\vbox{\offinterlineskip
\hrule
\halign{&\vrule#&\strut\quod\hfil#\quad\cr
height2pt&\omit&&\omit&&\omit&\cr
& $g$ && $N$ &&  Interval  &\cr
\noalign{\hrule}
& $41$ && $72$ && $[65,78]$  &\cr
} \hrule}
}}
\end{minipage}

\noindent\begin{minipage}{\textwidth}
\vspace{.2in}\noindent{\large
New entries for the tables over $\mathbf{F}_5$.}\\[.05in]
\smallskip
\vbox{
\bigskip\centerline{\def\quad{\hskip 0.6em\relax}
\def\quod{\hskip 0.5em\relax }
\vbox{\offinterlineskip
\hrule
\halign{&\vrule#&\strut\quod\hfil#\quad\cr
height2pt&\omit&&\omit&&\omit&\cr
& $g$ && $N$ &&  Interval  &\cr
\noalign{\hrule}
& $8$   && $24$ &&  $[22,28]$  &\cr
& $10$ && $31$ &&  $[27,33]$   &\cr
& $12$ && $36$ &&  $[33,38]$   &\cr
& $26$ && $60$ &&  $[\quad,68]$  &\cr
& $35$ && $72$ &&  $[68,85]$   &\cr
& $37$ && $80$ &&  $[72,89]$   &\cr
& $40$ && $72$ &&  $[\quad,94]$  &\cr
& $45$ && $96$ &&  $[88,104]$   &\cr
& $46$ && $81$ &&  $[75,106]$   &\cr
} \hrule}
}}
 (For genera where no good curve was known the lower bound on the interval is left empty).
\end{minipage}

\section{Details for the construction}\label{details}
In this section we give all data that are needed to construct the curves listed in Section \ref{results}. We begin by giving some of these constructions as detailed examples; then follows a list which in each case give just the the necessary details.

To make it easy for the reader to check the result, below we use Magma's notation for representing function fields and places; the field $F$ is represented as a finite separable extension of a rational function field $k(x)$. Every place of $F$ then corresponds either to a prime ideal in the integral closure of $k[x]$ or one in the integral closure of the valuation ring of the degree valuation.

\begin{example}[Genus 17 over $\FF_2$]
Let $F$ be the hyperelliptic genus 2 field defined by $y^2 + (x^2 + x + 1)y + x^5 + x^4 + x^2 + x$ over the rational field $\FF_2(x)$. Let $S=\{ (1/x, y/x^3), (x + 1, y + 1) \}$ and let $P=(x + 1, y + x + 1)$. We get a sequence of fields $F\subset F_S^{2P}\subset  F_S^{3P}\subset  F_S^{5P} $.

First, $\Cl_{2P}$ is isomorphic to $\ZZ/28\oplus\ZZ$. In it, $S$ generates a subgroup of index $4$. The field $F_S^{2P}$ has genus 7 and 10 rational places, the best possible. (This was already known). We get an explicit equation for this extension as
$T^4 + (x^3 +x)T^2 + (x^4 + 1)T + (x^6 + x^3 + x^2 + x)y + x^{16} + x^{12} + x^{11}
    + x^{10} + x^9 + x^8 + x^7 + x^5 + x^4 + 1$.

Next, $\Cl_{3P}$ is isomorphic to $\ZZ/56\oplus\ZZ$. The subgroup generated by $S$  in $\Cl_{3P}$ has index 8. The ray class field $F_S^{3P}$ has genus 17 and 18 rational places. This is a new record and the best possible, showing that $\NN_2(17)=18$. We can also get an explicit equation for $F_S^{3P}$  (which is very long, therefore given in the last section). Then we can construct the field without using class field theory, and then use Magma to compute its genus and number of rational places. This gives an independent verification that it really has genus 17 and 18 rational places.

Furthermore $\Cl_{5P}$ is isomorphic to $\ZZ/2\oplus \ZZ/112 \oplus \ZZ$. The field $F_S^{5P}$ has genus 45 and 34 rational places, the old record was 33 (in the next example we give one with $36$ places).
\end{example}

\begin{example}[Genus 45 over $\FF_2$]
Take the genus two field given as an extension of $\FF_2(x)$ by
$y^2 + (x^3 + x + 1)y + x^6 + x^5 + x^4 + x^2$. It has 4 rational places and 3 places of degree 2; use two of the latter to construct the divisor
$D=(x^2 + x + 1, y + x + 1) + 3(x^2 + x + 1, y + x^2 + x)$, and let $U$ be generated by the images of the three rational places
$\{ (x, y + x), (x + 1, y), (x + 1, y + 1) \}$ in $\Cl_D$. The index of $U$ is 12 and $F_U^D$ has genus 45 and 36 rational places, showing that $\NN_2(45)$ equals either $36$ or $37$.
\end{example}

\begin{example}[Genus 22 and 46 over $\FF_3$]
Let $F$ be given by $y^2 + 2x^6 + x^5 + 2x^4 + x^3 + 2x^2 + x + 2$, let $D=2(1/x, y/x^3 + 1) + 2(x, y + 2)$ and let
$S=\{ (x + 1, y + x + 2), (x + 2, y + x), (x + 2, y + x + 1)\}$. The field $F_S^D$ has genus 46 and 60 rational places. It has a subfield $K$ of genus 22 with 33 rational places. By listing subgroups of $\Cl_D$ this is easy to find; however, if we want to express $K$ as $F_T^D$ for some set of places $T$, then we have to use places of rather high degree: The easiest way to obtain such $T$ is to add the degree $5$ place $(x^5 + x^3 + x + 1, y + 2x^4 + x^3 + 2x)$ to $S$. Then $K= F_T^D\subset F_S^D$.

Explicitly, $F_T^D$ can be given as an extension of $F$ by the equations
\begin{align*}
 &T_1^3 + 2T_1 + (x^3 + 2x^2 + x + 1 + 1/x)y + x^6 + 2x^2 + 2x + 1/x  \\
 &T_2^3 + 2T_2 + (x^3 + 2x^2 + x + 2/x)y + x^6 + 2x^3 + 2x^2 + 2 + 2/x
\end{align*}
which we also use to verify that it has the claimed genus and number of places.
\end{example}

\subsection*{List of data sufficient to construct each curve}
Below we give sufficient details to construct each of the records that we found: For each entry in this list, there is first a finite field $\FF_q$ and a pair $(g,N)$; then a global function field $F/\FF_q$; then a divisor $D$ of this function field; then a set $S$ of places of $F$. The ray class field $F^D_S/\FF_q$  then has genus $g$ and $N$ rational places. This can be proved using the algorithms for computing ray divisor class groups which are given in \cite{Hess} and implemented for example in Magma. Also, in cases where it seems reasonable with respect to space we have included defining equation for the extension $F^D_S/F$, which we then also have used to give an independent verification that the field really has the claimed genus and number of rational places.

Although Magma is a proprietary software, there is a calculator \\  \emph{magma.maths.usyd.edu.au/calc/} that is free to use and sufficient to construct the records from the details given below.
\subsection*{Details for the curves over $\FF_2$}
\begin{itemize}
\item \underline{$\FF_2$: \qquad $(g,N)=(17, 18)$}
\begin{align*}
& F: y^2 + (x^2 + x + 1)y + x^5 + x^4 + x^2 + x \\
& D=3(z + 1, y + z + 1)  \\
& S=\{(1/z, y/z^3), (z + 1, y + 1) \}.
\end{align*}
Defining polynomial:
\begin{align*}
  T^8 +& (x^5 + x^4 + x + 1)T^6 + (x^7 + x^6 + x^5 + x^4 + x^3 + x^2 + x + 1)T^5  \\
    +& ((x^{17} + x^{15} + x^{14} + x^{13} + x^{12} + x^{10} + x^7 + x^5 + x + 1)y + (x^{22}\\
 +&  x^{20} + x^{19} + x^{14} + x^{12} + x^{11} + x^{10} + x^7 + x^5 + x^3 + x^2 + x))T^4  \\
    +&(x^{11} + x^{10} + x^9 + x^8 + x^3 + x^2 + x + 1)T^3  \\
+&((x^{31} + x^{29} + x^{28} +  x^{27} + x^{24} + x^{23} + x^{21}  \\
 +& x^{16} + x^{13} + x^{10} + x^9 + x^5 + x^4 + x^3 + x^2 + x)y  \\
 +&(x^{35} + x^{32} + x^{31} + x^{30} + x^{29} + x^{27} + x^{26} + x^{25} + x^{24} \\
  +&  x^{19} + x^{15} + x^{12} + x^{10} + x^6 + x^5 + x^4 + x^3 + x))T^2 \\
 +& ((x^{33} + x^{30} + x^{28} + x^{27} + x^{26} + x^{25} + x^{24} + x^{23}\\
 +& x^{21} + x^{20} + x^{19} + x^{18} + x^{17} + x^{15} + x^{14} + x^{13}  \\
 +&x^{12} + x^{11} + x^9 + x^4 + x^3 + 1)y +  (x^{37} + x^{35} + x^{34}  \\
+& x^{33} + x^{30} + x^{24} + x^{23} + x^{16} + x^{15} + x^{13} + x^6 + x^3 + x^2 + x))T \\
 +& (x^{54} + x^{53} + x^{51} + x^{50} + x^{47} + x^{45} + x^{44} + x^{43} + x^{42}\\
 +& x^{41} + x^{39} + x^{38} + x^{37} + x^{36} + x^{29} + x^{24} + x^{22} + x^{19} + x^{17} + x^{16}  \\
+& x^{15}+ x^{11} + x^{10} + x^7 + x^6 + x^2)y + x^{60} + x^{58} + x^{57} + x^{55} + x^{54}  \\
+& x^{53}+ x^{52} + x^{51} + x^{47} + x^{46} + x^{45} + x^{40} + x^{39} + x^{37}+ x^{35} + x^{33}  \\
+& x^{32} + x^{30} + x^{29} + x^{28} + x^{26} + x^{21} + x^{19} + x^{18} + x^{17} + x^{16} + x^{15}  \\
+&x^{13}+ x^{11} + x^9 + x^6 + x^4 + x^2 + 1.
\end{align*}

\item \underline{$\FF_2$: \qquad $(g,N)=(45,36)$}
\begin{align*}
&F:y^2 + (x^3 + x + 1)y + x^6 + x^5 + x^4 + x^2 \\
&D=(x^2 + x + 1, y + x + 1) + 3(x^2 + x + 1, y + x^2 + x)\\
&S=\{(x, y + x), (x + 1, y), (x + 1, y + 1) \}.
\end{align*}

\item \underline{$\FF_2$: \qquad $(g,N)=(46,36)$}
\begin{align*}
&F: y^2 + xy + x^5 + x^3 + x^2 + x \\
&D= (x^3 + x + 1) + (x^2 + x + 1, y + x + 1) + (x^2 + x + 1, y + 1) \\
&S=\{ (1/x, y/x^3), (x, y), (x + 1, y), (x + 1, y + x) \}.
\end{align*}

\item \underline{$\FF_2$: \qquad $(g,N)=(48,35)$}
\begin{align*} 
&F: y^2 + xy + x^5 + x \\
&D= (x^4 + x + 1, y + x^3 + x^2) + (x^4 + x + 1, y + x^3 + x^2 + x) \\
& \qquad+ 2(1/x,y/x^3) \\
&S=\{ (x, y), (x + 1, y + x + 1), (x + 1, y + 1), (x^2 + x + 1) \}.
\end{align*}

\end{itemize}

\subsection*{Details for the curves over $\FF_3$}
\begin{itemize}
\item \underline{$\FF_3$: \qquad $(g,N)=(17, 28)$}
\begin{align*}
&F: y^2 + x^5 + x^4 + x^2 + 2x \\
&D= (1/x, y/x^3) + (x + 1, y + 1) + 2(x^2 + 1, y)  \\
&S= \{(x + 1, y + 2), (x + 2, y + 1), (x + 2, y + 2) \}.
\end{align*}
Defining polynomials:
\begin{align*}
&T_1^4 + (x^3 + x)T_1^2 + 2x^6 + 2x^5 + x^4 + 2x^3 + 2x^2 \\
&T_2^2 + (x + 1)y + x^3 + x^2 + x + 1
\end{align*}

\item \underline{$\FF_3$: \qquad $(g,N)= (22, 33)$}
\begin{align*}
F&: y^2 + 2x^6 + x^5 + 2x^4 + x^3 + 2x^2 + x + 2 \\
D&=2(1/x, y/x^3+ 1) + 2(x, y + 2) \\
S&= \{(x + 1, y + x + 2), (x + 2, y + x), (x + 2, y + x + 1), \\
& (x^5 + x^3 + x + 1, y + 2x^4 + x^3 + 2x) \}.
\end{align*}
Defining polynomials:
\begin{align*}
& T_1^3 + 2T_1 + (x^3 + 2x^2 + x + 1 + 1/x)y + x^6 + 2x^2 + 2x + 1/x  \\
& T_2^3 + 2T_2 + (x^4 + 2x^3 + x^2 + 2)y/x + (x^7 + 2x^4 + 2x^3 + 2x + 2)/x
\end{align*}

\item \underline{$\FF_3$: \qquad $(g,N)=(33,48)$}
\begin{align*}
F&:y^2 + 2x^6 + x^5 + 2x^4 + x \\
D&=(x^2 + 1, y + x + 2) + (x^2 + 1, y + 2x + 1) \\
S&=\{ (1/x, y/x^3 + 1), (1/x, y/x^3 + 2), (x, y) \}.
\end{align*}

\item \underline{$\FF_3$: \qquad $(g,N)=(46,60)$}
\begin{align*}
&F:y^2 + 2x^6 + x^5 + 2x^4 + x^3 + 2x^2 + x + 2 \\
&D=2(1/x, y/x^3 + 1) + 2(x, y + 2) \\
&S=\{ (x + 1, y + x + 2), (x + 2, y + x), (x + 2, y + x + 1) \}.
\end{align*}

\end{itemize}

\subsection*{Details for the curves over $\FF_4$}
Below, $a$ is a primitive element of $\FF_4$.
\begin{itemize}
\item \underline{$\FF_4$: \qquad $(g,N)=(41,72)$}
\begin{align*}
&F:y^2 + (x^2 + x)y + x^5 + x^3 + a^2x^2 + a^2x \\
&D=(x + a, y + x) + (x + a, y + a^2) \\
&S=\{ (1/x, y/x^3), (x, y), (x + 1, y) \}.
\end{align*}
\end{itemize}

\subsection*{Details for the curves over $\FF_5$}
\begin{itemize}
\item \underline{$\FF_5$: \qquad $(g,N)=(8,24)$}
\begin{align*}
&F: y^2 + 4z^6 + 2z^5 + 2z^3 + z^2 + 2z \\
&D=(z^2 + 2z + 3, y + 4z + 4) + (z^2 + 4z + 2, y + 3z + 2) \\
&S=\{ (1/z, y/z^3 + 1), (1/z, y/z^3 + 4), (z, y), (z + 3, y + 4), \\
&\qquad (z + 1, y + 4), (z + 2, y + 2), (z + 4, y + 2), (z + 4, y + 3) \}.
\end{align*}
Defining polynomial:
\begin{align*}
T^3& + ((4z + 3)y + (3z^4 + 3z^2 + 2))T  \\
 &+(4z^3 + 3z^2 + 4z + 2)y + 4z^6 + z^4 + z^3 + z^2 + 2z + 2.
\end{align*}

\item \underline{$\FF_5$: \qquad $(g,N)=(10,31)$}
\begin{align*}
F&:y^2 + 4z^6 + 2z^5 + 3z^3 + 4z^2 + 1 \\
D&=3(z + 3, y + z) \\
S&=\{ (1/z, 1/z^3y + 1), (1/z, y/z^3 + 4), (z, y + z + 2), \\
&\qquad (z + 3, y + z + 1), (z + 1, y + 1), (z + 4, y + 4) \}.
\end{align*}
Defining polynomial:
\begin{align*}
T^5 +& 4T + (4z^3 + z^2 + 3z + 4)y/(z + 3) + \\
 +&(z^6 + 3z^5 + 4z^2 + 2z + 3)/(z + 3)
\end{align*}

\item \underline{$\FF_5$: \qquad $(g,N)=(12,36)$}
\begin{align*}
F&:y^2 + 3z^6 + z^5 + 2z^4 + 4z^3 + 4z^2 + 3z + 4 \\
D&=2(1/z) \\
S&=\{ (z, y + 1), (z, y + 4), (z + 2, y + z + 3), \\
& \quad (z + 2, y + z + 1), (z + 4, y + 2), (z + 4, y + 3) \}.
\end{align*}
Defining polynomials:
\begin{align*}
&x^2 + 4z^3 + 2z^2 + 4z + 1 \\
&w^3 + (3z^3 + z^2 + 2z + 1)w + z^5 + 3z^3 + 3z
\end{align*}

\item \underline{$\FF_5$: \qquad $(g,N)=(26,60)$}
\begin{align*}
F&: y^2 + z^6 + z^5 + 4z^4 + 4z^3 + 4z^2 + z + 1 \\
D&=(z^2 + z + 1, y + 3) + (z^2 + 3, y + 3z + 2)  \\
S&=\{(1/z, y/z^3 + 2), (z, y + 2), (z + 3, y + z), \\
& \quad (z + 2, y + 4) \}.
\end{align*}

\item \underline{$\FF_5$: \qquad $(g,N)=(35,72)$}
\begin{align*}
F&: y^2 + z^6 + 4z^5 + 2z^4 + 2z^2 + 4z + 1 \\
D&= (z^2 + z + 1, y + 2) + (z + 1)\\
S&=\{(1/z, y/z^3 + 2), (z, y + 2), (z + 3, y + 2),  \\
& \quad (z + 2, y + z + 1), (z + 4, y + 1), (z + 4, y + 4)  \}.
\end{align*}

\item \underline{$\FF_5$: \qquad $(g,N)=(37,80)$}
\begin{align*}
F&: y^2 + z^5 + z^4 + 2z^3 + z^2 + 4z \\
D&=3(z, y) \\
S&=\{ (1/z, y/z^3), (z + 1, y), (z + 4, y + 1),   \\
& \quad (z + 4, y + 4)  \}.
\end{align*}

\item \underline{$\FF_5$: \qquad $(g,N)=(40,72)$}
\begin{align*}
F&:y^2 + 2z^5 + z^4 + 2 \\
D&=(z + 3) + (z)\\
S&=\{(1/z, y/z^3), (z + 1, y + 2), (z + 1, y + 3),    \\
& \quad  (z + 4, y) \}.
\end{align*}

\item \underline{$\FF_5$: \qquad $(g,N)=(45, 96)$}
\begin{align*}
F&:y^2 + 2z^6 + 4z^4 + 3z^2 + 1 \\
D&=2(1/z)\\
S&=\{ (z + 3, y), (z + 1, y), (z + 2, y), (z + 4, y) \}.
\end{align*}

\item \underline{$\FF_5$: \qquad $(g,N)=(46, 81)$}
\begin{align*}
F&:y^2 + z^6 + 2z^5 + 2z^4 + z^3 + 2z^2 + 2z + 1 \\
D&=2(z^2 + 4z + 2, y + z^2 + 4)\\
S&=\{ (z, y + 2), (z + 3, y + 3), (z + 1, y + 3) \}.
\end{align*}

\end{itemize}

\begin{acknowledgements}
The author thanks the Wenner-Gren Foundations for
financial support; the Korteweg-de Vries Institute at the University of
Amsterdam for providing office space; and Gerard van der Geer for his
hospitality and for helpful conversations and comments. 
Thanks are also due to Jan-Erik Roos at Stockholm University for giving me time on his computer to run Magma, and to Samuel Lundqvist at Stockholm University for help with this.
\end{acknowledgements}

\comment{\bibliographystyle{amsplain}
\bibliography{biblio}}

\providecommand{\bysame}{\leavevmode\hbox
to3em{\hrulefill}\thinspace}

\end{document}